\documentclass[leqno,11pt]{amsart}
\usepackage{amsmath}
\usepackage{amssymb,amscd}
\usepackage{amsthm}
\usepackage{latexsym}
\usepackage{hyperref}
\usepackage[myheadings]{fullpage}
\usepackage{color}

\newtheorem{theorem}{Theorem}[section]
\newtheorem{lemma}[theorem]{Lemma}
\newtheorem{proposition}[theorem]{Proposition}
\newtheorem{corollary}[theorem]{Corollary}
\newtheorem{remark}[theorem]{Remark}
\newtheorem{example}[theorem]{Example}

\newtheorem{definition}[theorem]{Definition}
\newtheorem{question1}{Question}

\def\opn#1#2{\def#1{\operatorname{#2}}} 

\opn\PF{PF}
\opn\F{F}
\opn\G{G}
\opn\RF{RF}

\title{Nearly Gorenstein numerical semigroups with five generators have bounded type}

\author{Alessio Moscariello}
\address{A. Moscariello - Dipartimento di Matematica e Informatica -  Universit\`a di Catania - Viale Andrea Doria 6 - 95125 Catania - Italy.}
\email{alemoscariello@hotmail.it}

\author{Francesco Strazzanti}
\address{F. Strazzanti - Dipartimento di Scienze Matematiche e Informatiche, Scienze Fisiche e Scienze della Terra - Universit\`a di Messina - Viale Ferdinando Stagno d’Alcontres 31 - 98166 - Messina - Italy.}
\email{francesco.strazzanti@gmail.com}

\thanks{The first author was funded by the project “Propriet\`a locali e globali di anelli e di variet\`a algebriche”-PIACERI 2020–22, Universit\`a degli Studi di Catania.
The second author was supported by the "National
Group for Algebraic and Geometric Structures, and their Applications" (GNSAGA - INdAM)}


\begin{document}

\begin{abstract}
We prove that the type of nearly Gorenstein numerical semigroups minimally generated by $5$ integers is bounded. In particular, if such a semigroup is not almost symmetric, then its type is at most $40$. Finally, we make some considerations in higher embedding dimension.
\end{abstract}


\maketitle

\section*{Introduction}

Numerical semigroups $S$ are additive submonoids of $\mathbb{N}$ for which $\mathbb{N}\setminus S$ is finite. Despite their simple structure, they are source of interesting and challenging problems; indeed, they appear in several mathematical areas, like Commutative Algebra, Algebraic Geometry, Number Theory or Coding Theory, and many problems arise directly from these fields. For instance, it is possible to associate a numerical semigroup with a one-dimensional analytically irreducible ring and their property are closely related. The simplest example of analytically irreducible ring is the numerical semigroup ring $\mathbb K[[S]]=\mathbb{K}[[t^s \mid s \in S]]$ associated with the numerical semigroup $S$ and a field $\mathbb K$. Many invariants of $\mathbb K[[S]]$ can be easily read off from the arithmetic of $S$, like embedding dimension, multiplicity, Hilbert function, Gorenstein property or Cohen-Macaulay type. For example the embedding dimension is the minimal number of generators of $S$, while the Cohen-Macaulay type is equal to the number of the so-called pseudo-Frobenius numbers of the semigroup $S$, and we will denote it by $t(S)$, since does not depend on the field $\mathbb{K}$. It is easy to see that a numerical semigroup has at least two generators and in this case $t(S)=1$ because $\mathbb K[[S]]$ is a complete intersection. When $S$ has three generators, a celebrated result by Herzog \cite{H} ensures that $t(S) \leq 2$. Nevertheless, an example due to Backelin \cite{FGH} (see also Example \ref{Backelin}) shows that $t(S)$ is not bounded if $S$ has four generators. Numata \cite{Nu} conjectured that if we assume that a four generated numerical semigroup $S$ is almost symmetric, then $t(S) \leq 3$. This conjecture was proved to be true by the first author \cite{M}, who also proved that a five generated almost symmetric numerical semigroup has type at most $473$ \cite{M5}. Even if this bound seems large, and it is probably not sharp, it means that a finite bound exists. Indeed, in general it is not known if the type of an almost symmetric semigroup can be bounded by a function involving only the minimal number of generators of $S$; the problem is widely open even in the 6-generated case. Almost symmetric semigroups were introduced by Barucci and Fr\"oberg in \cite{BF} as a generalization of symmetric semigroups. Recall that $S$ is symmetric if and only if $\mathbb K[[S]]$ is Gorenstein, while $S$ is almost symmetric precisely when $\mathbb K[[S]]$ is an almost Gorenstein ring. In the one-dimensional case, a larger class is given by nearly Gorenstein ring, which were introduced and systematically studied by Stamate, Hibi, and Herzog in \cite{HHS}, even if this notion had already appeared in literature \cite{D,HV}. If $\mathbb K[[S]]$ is nearly Gorenstein, we also say that $S$ is nearly Gorenstein; see for instance \cite{HHS2,MS} for several properties of these semigroups. Note that an almost symmetric semigroup is also nearly Gorenstein. In \cite{MS} it is proved that the type of a $4$-generated nearly Gorenstein numerical semigroup is at most $3$. When there are more minimal generators, there are the following open questions:

\begin{question1} \cite[Question 3.7]{S}
Is there any bound for the type of $\mathbb K[[S]]$ in terms of the embedding dimension of $S$ when $S$ is nearly Gorenstein?
\end{question1}

\begin{question1} \cite[Question 4.2]{MS}
Let $S$ be a nearly Gorenstein numerical semigroup with five minimal generators. Is it true that $t(S) \leq 5$ and that the equality is attained only if $S$ is almost symmetric?
\end{question1}

In this paper we focus on the nearly Gorenstein numerical semigroups with five generators that are not almost symmetric, and prove that their type is indeed bounded by showing that it is at most $40$. The main tool are the notions of NG-vector and row factorization matrices. Moreover, we provide some ideas on how to use row factorization matrices to construct examples with six generators and type not bounded, even if it is still not clear how to produce almost symmetric or nearly Gorenstein numerical semigroups in this way.  

\medskip

The structure of the paper is the following. In Section 1, we recall some definitions and known results about NG-vectors and row factorization matrices. In Section 2, we prove the bound for the type of $5$-generated nearly Gorenstein numerical semigroups, see Theorem \ref{main}. In Theorem \ref{3distinct}, we also prove that for a particular family of nearly Gorenstein numerical semigroups the upper bound for their type is $5$. Finally, in the last section we give examples and some ideas for constructing numerical semigroups with fixed embedding dimension and type not bounded. In Remark \ref{Numerical duplication}, we also show that the type of almost symmetric numerical semigroups with $\nu$ generators is not bounded by $2\nu+a$ with $a \in \mathbb{Z}$.

\section{Preliminaries}

Let $S$ be a numerical semigroup, that is an additive submonoid of $\mathbb N$ such that $\mathbb{N}\setminus S$ is finite. We write $S=\langle n_1, \dots, n_\nu \rangle$ to mean that $\{n_1, \dots, n_\nu\}$ is the unique minimal system of generators of $S$, i.e., $S=\{a_1 n_1+\dots+a_\nu n_\nu \mid a_i \in \mathbb{N} \text{ for every } i=1, \dots, \nu\}$ and it is not possible to remove some $n_i$ in the list obtaining $S$. The integer $\nu$, that is the minimal number of generators of $S$, is called {\it embedding dimension} of $S$.
We denote by $\F(S)$ the {\it Frobenius number} of $S$, i.e., the biggest integer that is not in $S$.
An integer $f \notin S$ is said to be a {\it pseudo-Frobenius number} of $S$ if $f+s \in S$ for every $s \in M(S)$, where $M(S)=S \setminus \{0\}$. We denote the set of pseudo-Frobenius numbers of $S$ by $\PF(S)$; its cardinality $t(S)$ is the {\it type} of $S$. We introduce the order $\leq_S$ in $\mathbb{Z} \setminus S$ by setting $x \leq_S y$ if there exists $s \in S$ such that $y=x+s$. Note that the pseudo-Frobenius of $S$ are exactly the maximal elements of $\mathbb{Z}\setminus S$ with respect to $\leq_S$, thus if $f,f' \in \PF(S)$ and $f \leq_S f'$, it follows that $f=f'$.
We refer to \cite{ADG} or \cite{RG} for more details about numerical semigroups.

In this paper we are interested in the type of some particular class of semigroups that generalizes symmetric numerical semigroups. Recall that the {\it canonical ideal} of $S$ is $K(S)=\{s \in \mathbb{N} \mid \F(S)-s \notin S\}$, which contains $S$. The semigroup $S$ is said to be symmetric when $S=K(S)$, whereas it is called almost symmetric if $K(S)+M(S) \subseteq M(S)$. The central notion in this paper is the nearly Gorenstein property, which was defined in \cite{HHS} in a more general setting. In our context we may say that $S$ is {\it nearly Gorenstein} if $K(S)+(S-K(S)) \supseteq M(S)$, where $S-K(S)=\{x \in \mathbb{Z} \mid x+K(S) \subseteq S\}$. However, this definition are not effective for our purpose, so we will always use the characterizations in the next proposition, where the first one is essentially a reformulation of part of \cite[Theorem 2.4]{N}. 

\begin{proposition}\cite[Proposition 1.1]{MS} \label{Nearly Gorenstein}
The following statements hold:
\begin{enumerate}
\item $S$ is almost symmetric if and only if $n_i+\F(S)-f \in S$ for all $f \in \PF(S)$ and all $i=1,\dots, \nu$;
\item $S$ is nearly Gorenstein if and only if for every $i=1, \dots, \nu$ there exists $f_i \in \PF(S)$ such that $n_i+f_i-f \in S$ for all $f \in \PF(S)$. 
\end{enumerate}
In particular, an almost symmetric numerical semigroup is nearly Gorenstein.
\end{proposition}

The characterization of the nearly Gorenstein property given in the previous proposition led to introduce the following notion.

\begin{definition} \rm
Let $S=\langle n_1, \dots, n_{\nu}\rangle$, where $n_1 < \dots < n_{\nu}$ are minimal generators. We call a vector $(f_1,\dots, f_{\nu}) \in \PF(S)^{\nu}$ {\it nearly Gorenstein vector} for $S$, briefly {\rm NG}-vector, if $n_i+f_i-f \in S$ for every $f \in \PF(S)$ and every $i=1, \dots, \nu$.
\end{definition}

By Proposition \ref{Nearly Gorenstein}, the existence of an NG-vector is equivalent to the nearly Gorenstein property of $S$, whereas $S$ is almost symmetric if and only if it admits the NG-vector $(\F(S), \dots, \F(S))$. Notice that a numerical semigroup may have more NG-vectors.
These vectors were introduced in \cite{MS}, where it is possible to find several properties and examples. In the next proposition we recall two results proved in \cite{MS} and prove a new one.

\begin{proposition} \label{minimum index}
Let $(f_1, \dots, f_{\nu})$ be an {\rm NG}-vector for $S$. The following hold:
\begin{enumerate}
\item $f_1=\F(S)$.
\item If $h$ is the smallest index for which $f_h \neq \F(S)$, then $f_h=\F(S)-n_h+n_\ell$ for some $\ell < h$.
\item Assume there is another index $h'\neq h$ such that $f_{h'} \neq \F(S)$ and let $h'$ the minimum of such indices. Then, either $f_{h'}=\F(S)-n_{h'}+n_{\ell'}$ for some $\ell'<h'$ or $f_{h'}=\F(S)-n_{h'}+a_h n_h$ for some positive integer $a_h$.
\end{enumerate}
\end{proposition}

\begin{proof}
(1) and (2) are proved in \cite[Propostion 1.3]{MS}, while we prove (3) here. Since $f_{h'}<\F(S)$, we have $n_{h'}+f_{h'}-\F(S)<n_{h'}$, and then $n_{h'}+f_{h'}-\F(S)= a_1n_1+\dots+a_{h'-1}n_{h'-1}$ for some $a_1,\dots,a_{h'-1} \in \mathbb{N}$ with at least one of them positive.

Assume first that $a_{\ell'} >0$ for some $\ell'\neq h$. Then, $$n_{h'} = \F(S)-f_{h'}+n_{\ell'}+ a_1n_1 + \cdots + (a_{\ell'}-1)n_{\ell'} + \dots +a_{h'-1}n_{h'-1}.$$ Since $\F(S)=f_{\ell'}$, it follows that $\F(S)-f_{h'}+n_{\ell'} \in S \setminus \{0\}$. On the other hand, $n_{h'}$ is a minimal generator, and then $a_{\ell'}=1$ and $a_{i}=0$ if $i \neq \ell'$, i.e., $f_{h'}=\F(S)-n_{h'}+n_{\ell'}$.

Finally, if $a_{\ell'}=0$ for every $\ell'\neq h$, then $f_{h'}=\F(S)-n_{h'}+a_h n_h$.
\end{proof}

We also recall other two results from \cite{MS} that we will use later.

\begin{proposition} \cite[Proposition 2.8]{MS} \label{f1-fj}
Let $(f_1, \dots, f_{\nu})$ be an {\rm NG}-vector for $S$ and suppose that $f_1, f_2, \dots, f_i$ are pairwise distinct for some $i \leq \nu$. The following statements hold:
\begin{enumerate}
	\item If $f \in \PF(S) \setminus \{f_1, \dots, f_i\}$, then $f_1-f=-n_1+a_{i+1}n_{i+1}+ \dots + a_{\nu}n_{\nu}$, with $a_j \in \mathbb{N}$ for every $j=i+1, \dots, \nu$;
    \item $f_1-f_{j} =n_{j}-n_1$ for every $j=1, \dots, i$.
\end{enumerate}
\end{proposition}

\begin{corollary} \cite[Corollary 2.9]{MS} \label{Distinct}
If $(f_1, \dots, f_{\nu})$ is an {\rm NG}-vector for $S$, then there exist at least two different indices $i$ and $j$ such that $f_i=f_j$. Moreover, if $f_1, \dots, f_{\nu-1}$ are pairwise distinct, then $\PF(S)=\{f_1, \dots, f_{\nu-1}\}$ and $t(S)=\nu-1$.
\end{corollary}

\subsection{RF-matrices}

Let $S=\langle n_1, \dots, n_\nu \rangle$. The notion of $\RF$-matrix was defined in \cite{M} and has proved to be a powerful tool in investigating the type and the defining ideal of a numerical semigroup, see for instance \cite{E,HW,M5}. The idea is that an integer $f \notin S$ is a pseudo-Frobenius number if and only if $f+n_i \in S$ for every $i$; this means that $f+n_i=\sum_{j=1}^\nu a_{ij}n_j$, where $a_{ii}=0$ (otherwise $f \in S$), and we can put all these factorizations in a matrix.

\begin{definition} \rm
We say that a square matrix $A=(a_{ij})$ of order $\nu$ is a {\it row factorization matrix} for $f \in \PF(S)$, briefly $\RF^+$-matrix, if $a_{ii}=-1$, $a_{ij} \in \mathbb{N}$ when $i \neq j$, and $f=\sum_{j=1}^{\nu}a_{ij}n_j$ for every $i=1, \dots, \nu$.    
\end{definition}

Following \cite{MS}, in the previous definition we have written $\RF^+$, instead of $\RF$ as in \cite{M}, to avoid confusion with another matrix that we are going to introduce.

Assume that $S$ is nearly Gorenstein, and let $(f_1, \dots, f_{\nu})$ be an NG-vector for $S$. If $f \neq f_i$ for every $i$, we have
\[
n_i+f_i-f=\sum_{j=1}^{\nu}b_{ij}n_j 
\]
with $b_{ij} \geq 0$ and $b_{ii}=0$, since $f_i-f \notin S$. Thus, we can define another matrix.

\begin{definition} \rm
Let $(f_1, \dots, f_{\nu})$ be an {\rm NG}-vector for $S$ and let $f \in \PF(S)\setminus \{f_1, \dots, f_\nu\}$. We say that a square matrix $B=(b_{ij})$ of order $\nu$ is an $\RF^-$--matrix for $f$ if $b_{ii}=-1$, $b_{ij} \in \mathbb{N}$ when $i \neq j$, and $f_i-f=\sum_{j=1}^{\nu}b_{ij}n_j$ for every $i=1, \dots, \nu$.  
\end{definition}

Note that, even if we fix an NG-vector, there could be more $\RF^+$ and $\RF^-$--matrices associated with a pseudo-Frobenius number $f$. However, these matrices are always related by the following fundamental result.

\begin{lemma} \cite[Lemma 2.3]{MS} \label{Coppie}
Let $S$ be nearly Gorenstein and let $f\in \PF(S)$. Also, let $A=(a_{ij})$ and $B=(b_{ij})$ be an $\RF^+$ and an $\RF^{-}$--matrix for $f$ respectively. Then, $a_{jk}b_{kj}=0$ for every $j\neq k$.  
\end{lemma}

\begin{remark} \label{first zero} \rm
Let $h$ and $\ell$ be as in Proposition \ref{minimum index} and let $f \in \PF(S)\setminus \{\F(S), f_h\}$. In every $\RF^-$--matrix for $f$, the entries with indices $(h,\ell)$ and $(\ell,h)$ are zero. Indeed, if $f_h-f=-n_h+ \sum_{j=1}^{h-1} b_{hj}n_j + \sum_{j=h+1}^{\nu} b_{hj}n_j$, then by Proposition \ref{minimum index}
\begin{align*}
\F(S)-f&=(\F(S)-f_h) + (f_h-f)=(-n_\ell+n_h)+\left(-n_h+\sum_{j=1}^{h-1} b_{hj}n_j + \sum_{j=h+1}^{\nu} b_{hj}n_j\right)= \\
&=(b_{h\ell}-1)n_\ell+ \sum_{j \in \{1,\dots, \nu\}\setminus \{h,\ell\}} b_{hj}n_j
\end{align*}
and this implies that $b_{h\ell}=0$, otherwise $\F(S)-f$ would be in $S$. Starting from a factorization of $\F(S)-f$ and proceeding in the same way it is straightforward to see that also the $(\ell, h)$-entry is zero. Moreover, this also means that there is a $\RF^-$--matrix of $f$ for which the $\ell$-row and $h$-row are exactly the same except for the entries in the $\ell$ and $h$ columns that are $0$ or $-1$.
\end{remark}

\section{On the type of 5-generated nearly Gorenstein semigroups}

In \cite{M5} it is proved that the type of $5$-generated almost symmetric numerical semigroups is bounded by $473$. In this section we will show a similar result for the nearly Gorenstein case.
Throughout the section we assume that $S$ is a nearly Gorenstein numerical semigroup minimally generated by the integers $n_1 <\ldots <n_5$, and that $S$ is not almost symmetric. We also fix an NG-vector $(f_1, f_2, \dots,f_5)$ for $S$.  
We start with some easy but crucial observations. 

\begin{remark} \rm \label{two types}
Let $f \in \PF(S)\setminus \{f_1, \dots, f_5\}$ and denote by $\RF^+(f)$ and $\RF^-(f)$ an $\RF^+$ and an $\RF^-$-- matrix for $f$ respectively. Note that in both $\RF^+(f)$ and $\RF^-(f)$ there are no rows with four zeroes.
    
$\bullet$ If a row of $\RF^+(f)$ contains three zeroes, then there exist two indices $i,j$ such that $f=\lambda n_j-n_i$ for some $\lambda \in \mathbb N$. We set $M_{i,j}=\lambda_{ij} n_j-n_i$, where $\lambda_{ij}=\max\{\lambda \mid \lambda n_j-n_i \notin S\}$. Since $f$ is maximal with respect to $\leq_S$, if $f=\lambda n_j-n_i$ for some positive $\lambda$, then $f=M_{i,j}$. In particular, there is at most one $f$ that can be written in this form (with the same indices).

$\bullet$ If a row of $\RF^-(f)$ contains three zeroes, then there exist two indices $i,j$ such that $f_i-f=\lambda n_j-n_i$ for some $\lambda \in \mathbb N$. Again there is at most one $f$ with this form because if $f_i-f'=\mu n_j-n_i$ for some $\mu \in \mathbb N$, then as before we get either $f \leq_S f'$ or $f' \leq_S f$, which implies $f=f'$.

$\bullet$ Assume that $f$ is not in the NG-vector of $S$. If in both $\RF^+(f)$ and $\RF^-(f)$ there are no rows with three zeroes, then in the two matrices there are at most 20 zeroes. On the other hand, Lemma \ref{Coppie} immediately implies that there are at least 20 zeroes, and then, using again Lemma \ref{Coppie}, we easily get that in every row and every column of  $\RF^+(f)$ and $\RF^-(f)$ there are exactly two zeroes.
\end{remark}

Recall that we have fixed an NG-vector $(f_1, \dots, f_5)$ for $S$. In light of the previous remark, in the following definition we divide the set $\PF(S)\setminus \{f_1, \dots, f_5\}$ in two sets. 

\begin{definition} \rm
Let $S$ be a nearly Gorenstein semigroup with emebedding dimension $\nu$ and NG-vector $(f_1, \dots, f_\nu)$. We denote by $\PF_1(S)$ the set of the pseudo-Frobenius numbers of $S$ different from $f_i$, for $i=1, \dots, \nu$, for which there exists at least a row of at least one $\RF^+$--matrix or $\RF^-$--matrix containing $\nu-2$ zeroes. We also denote by $\PF_2(S)$ the pseudo-Frobenius numbers of $S$  different from $f_i$, for $i=1, \dots, \nu$, that are not in $\PF_1(S)$.
\end{definition}

Clearly, by definition $t(S)=|\PF_1(S)|+|\PF_2(S)|+|\{f_1,f_2,f_3,f_4,f_\nu\}|$. In order to have a bound for $t(S)$ when $\nu=5$, in the following two subsections we study the cardinalities of the sets $\PF_1(S)$ and $\PF_2(S)$. It immediately follows from the definition that the cardinality of $\PF_1(S)$ is always bounded by $40$, whereas it is not clear whether there is a bound for $|\PF_2(S)|$ valid for all nearly Gorenstein numerical semigroups $S$ with $5$ minimal generators. Therefore, we start by studying the set $\PF_2(S)$ in the next subsection.

\subsection{A bound for $|\PF_2(S)|$} We begin by proving a lemma that will be crucial to show that $|\PF_2(S)|$ is indeed bounded.

\begin{lemma}\label{PF2}
Assume that $S=\langle n_1,\ldots,n_5 \rangle$ is nearly Gorenstein and let $f,f'\in \PF_2(S)$ be such that
	\begin{alignat}{2}
		f &= a_{ij}n_j+a_{il}n_l-n_i &= a_{pj}n_j+a_{pq}n_q-n_p\\ 
		f' &= b_{ij}n_j+b_{il}n_l-n_i &= b_{pj}n_j+b_{pq}n_q-n_p 
	\end{alignat}
with $\{i,j,l,p,q\} = \{1,\ldots,5\}$. Assume also that all the coefficients $a_{mn}$ and $b_{mn}$ are positive for every $m,n$ and $a_{ij} \ge a_{pj}$, $b_{ij} \ge b_{pj}$. Then $f=f'$.
\end{lemma}

\begin{proof}
Without loss of generality assume $a_{pq} \ge b_{pq}$. By (2) it follows that 
$$b_{pq}n_q=(b_{ij}-b_{pj})n_j+b_{il}n_l+n_p-n_i.$$
Using this equality in (1) we get
	\begin{alignat}{2} 
		f &=a_{pj}n_j+a_{pq}n_q-n_p \nonumber \\ &=a_{pj}n_j+(a_{pq}-b_{pq})n_q-n_p+(b_{ij}-b_{pj})n_j+b_{il}n_l+n_p-n_i \nonumber \\ &=(a_{pj}+b_{ij}-b_{pj})n_j+b_{il}n_l+(a_{pq}-b_{pq})n_q-n_i.
	\end{alignat}
Since $j \neq l$, $j \neq q$, $q \neq l$, $a_{pj}+b_{ij}-b_{pj} \ge a_{pj} > 0$, $b_{il} > 0$ and $a_{pq}-b_{pq} \ge 0$, the equality in (3) describe the $i$-th row of an $\RF^+$--matrix. By assumption, this row has to contain exactly $2$ positive coefficients, then $a_{pq}=b_{pq}$.
	
Therefore $f-f'=(a_{pj}-b_{pj})n_j$, and thus either $f \le_S f'$ or $f' \le_S f$, which imply $f=f'$.
\end{proof}

Now we can prove that the cardinality of $|\PF_2(S)|$ is at most 6.

\begin{proposition}
If $S=\langle n_1, \dots, n_5 \rangle$ is nearly Gorenstein, but not almost symmetric, then $|\PF_2(S)| \leq 6$.
\end{proposition}

\begin{proof} 
Let $f \in \PF_2(S)$, let $\RF^+(f)=(a_{ij})$ and $\RF^-(f)=(b_{ij})$ be an $\RF^+$ and $\RF^-$--matrix for $f$ respectively. Using the same notation as in Remark \ref{first zero}, the remark implies that $b_{h \ell}=b_{\ell h}=0$ and $b_{hj}=b_{\ell j}$ for every $j \not \in \{h,\ell\}$. 
In particular, there exist two indices $p,q$ (different from $h,\ell$) such that $b_{hp},b_{hq},b_{\ell p},b_{\ell q}$ are positive. Since $f \in \PF_2(S)$, it follows that $a_{ph}=a_{qh}=a_{p\ell}=a_{q\ell}=0$.

Now let $i$ be the index different from $p,q,h,\ell$. Then
\begin{equation} 
f=a_{pi}n_i+a_{pq}n_q-n_p=a_{qi}n_i+a_{qp}n_p-n_q.
\end{equation}
If there are three distinct pseudo-Frobenius numbers in $\PF_2(S)$ satisfying (4) with the same indices, then at least two of them satisfy the assumptions of Lemma \ref{PF2}; thus, two of them are equal and this yields a contradiction. This means that for every choice of the couple $(p,q)$ there are at most two elements in $\PF_2(S)$. On the other hand, since $p$ and $q$ are different from $h$ and $\ell$, there are three possible couples $(p,q)$. Therefore, $|\PF_2(S)| \le 6$.
\end{proof}

\subsection{Bounds for $|\PF_1(S)|$} Given $f \in \PF_1(S)$, we denote by $\RF^+(f)$ and $\RF^-(f)$ an $\RF^+$ and an $\RF^-$--matrix for $f$ respectively. Since $f \in \PF_1(S)$, we may assume that at least one row in either $\RF^+(f)$ or $\RF^-(f)$ contains three zeroes. This means that either $f=M_{i,j}=\lambda_{i,j} n_j-n_i$ or $f_i-f=\lambda n_j-n_i$ for some $j \neq i$ and $\lambda \in \mathbb{N}$. As noticed in Remark \ref{two types}, if we fixed the indices $i$ and $j$, there can be only one $f \in \PF(S)$ in the form $f=M_{i,j}=\lambda_{i,j} n_j-n_i$ and only one $f' \in \PF(S)$ in the form $f_i-f'=\lambda n_j-n_i$. This immediately implies that $|\PF_1(S)|\leq 40$. However, by Proposition \ref{minimum index} both $\F(S)-f_h$ and $f_h-\F(S)$ can be written in this way and so $|\PF_1(S)|\leq 38$ because $\F(S),f_h \notin \PF_1(S)$ by definition. Moreover, if in the NG-vector of $S$ there are at least two entries different from $\F(S)$, then Proposition \ref{minimum index}(3) implies that $|\PF_1(S)|\leq 37$. We now prove some results that will allow us to improve this bound.

\begin{remark} \rm
Note that $M_{i,k} \neq M_{j,k}$ when $i \neq j$. Indeed, if we assume that $\lambda_{ik}n_k-n_i=\lambda_{jk}n_k-n_j$ and $\lambda_{ik}\leq \lambda_{jk}$, then $n_j=(\lambda_{jk}-\lambda_{ik})n_k+n_i$. Since $n_j$ is a minimal generator, this implies that $\lambda_{jk}-\lambda_{ik}=0$ and $i=j$.  
\end{remark}

\begin{lemma}\label{same2}
Let $S=\langle n_1,\ldots,n_\nu \rangle$ be nearly Gorenstein. Assume that there exist $f,f' \in \PF_1(S)$ and three different indices $i,p,q$ for which $f=M_{p,i}=\lambda_{pi}n_i-n_p$ and $f'=M_{q,i}=\lambda_{qi}n_i-n_q$ with $\lambda_{pi} \ge \lambda_{qi}$. If $\RF^-(f)=(b_{ij})$, then $b_{qp}=0$.
\end{lemma}

\begin{proof}
Assume by contradiction that $b_{qp} \neq 0$. Looking at the $q$-th row of $\RF^-(f)$ we have 
$$f_q=f+(f_q-f)=\lambda_{pi}n_i-n_p + \sum_{j=1}^\nu b_{qj}n_j,$$
and, since $b_{qp} \ge 1$ and $\lambda_{pi} \ge \lambda_{qi}$, we get $f' \le_S f_q$, which is a contradiction.
\end{proof}

\begin{lemma}\label{same2+}
Let $S=\langle n_1,\ldots, n_5 \rangle$ be nearly Gorenstein. Suppose that there exist $s \in \{1,2,3,4,5\}$ and $\rho_1,\ldots,\rho_t \in \PF_1(S)$ such that for every $j$ there is an index $r_j$ for which $\rho_j=M_{r_j,s}=\lambda_{r_j,s} n_s - n_{r_j}$. Assume also that $\lambda_{r_1,s} \le \lambda_{r_2,s}\le \cdots \le \lambda_{r_t,s}$.
Then, for every $j=1,\ldots,t$, among the rows of $\RF^+(\rho_j)$ and $\RF^-(\rho_j)$ there are at least $j-1$ of them containing $3$ zeroes. 
\end{lemma}
  
\begin{proof}
Fix $j$, and consider $\RF^+(\rho_j)=(a_{pq})$ and $\RF^-(\rho_j)=(b_{pq})$. By applying Lemma \ref{same2} to $\rho_j$ and $\rho_i$ for every $i < j$, we get $b_{r_i,r_j}=0$ for every $i < j$. On the other hand, $\rho_j=\lambda_{r_j,s}n_s-n_{r_j}$ describe the $r_j$-th row of an $\RF^+$--matrix of $\rho_j$, thus we may assume that $a_{r_j,r_i}=0$ because $s \neq r_i$ by construction. 
Consider now the $20$ couples $\{a_{pq},b_{qp}\}$. By Lemma \ref{Coppie}, at least one element in every couple is zero; moreover, in the $j-1$ couples $\{a_{r_j,r_i},b_{r_i,r_j}\}$ with $i < j$ both elements are zero. This means that there are at least $20+j-1$ zeroes in $\RF^+(\rho_j)$ and $\RF^-(\rho_j)$. Since in every row there is at most three zeroes, among the rows of $\RF^+(\rho_j)$ and $\RF^-(\rho_j)$ there are at least $j-1$ of them containing $3$ zeroes. 
  \end{proof}
  
\begin{lemma}\label{mu}
Let $S=\langle n_1, \dots, n_5 \rangle$ be nearly Gorenstein but not almost symmetric. Define $\mu_s=|\{i \mid  M_{i,s} \in \PF_1(S)\}|$. Then
$$|\PF_1(S)| \le 38-\sum_{s=1}^5 \binom{\mu_s-1}{2}.$$
\end{lemma}

\begin{proof}
Fix an index $s$. There exist $\rho_1,\ldots,\rho_{\mu_s} \in \PF_1(S)$ satisfying the assumptions of Lemma \ref{same2+}. If  $t \le \mu_s$, among the rows of $\RF^+(\rho_t)$ and $\RF^-(\rho_t)$ there are at least $t-1$ of them containing three zeroes. To each of these $t-1$ rows corresponds  a writing of $\rho_t$ or $f_i-\rho_t$ as $\lambda n_j - n_i$, for some couple of indices $(i,j)$. In particular, if $t \ge 3$ there are $t-1$ couple $(i,j)$ for which either $\rho_t=\lambda n_j-n_i$ or $f_i-\rho_t=\lambda n_j - n_i$. Note that at least one element $M_{i,j}$ is associated with each $\rho_t$ by definition of $\rho_t$, but if $t \ge 3$ then there are at least $t-1$ of such writings associated with $\rho_t$. Since there are at most $38$ of such writings in $\PF_1(S)$ and every element of $\PF_1(S)$ is written in this way, we get 
\[
|\PF_1(S)| \leq 38-\sum_{s=1}^5 \sum_{t=3}^{\mu_s} (t-2)= 38-\sum_{s=1}^5 \binom{\mu_s-1}{2}.\qedhere \]
\end{proof}

\begin{proposition} \label{PF1}
Let $S=\langle n_1, \dots, n_5 \rangle$ be nearly Gorenstein but not almost symmetric. Then $|\PF_1(S)| \le 31$. Moreover, if in the {\rm NG}-vector of $S$ there are at least two entries different from $\F(S)$, then $|\PF_1(S)| \le 30$.
\end{proposition}

\begin{proof}
As in the previous lemma, let $\mu_s=|\{i \mid  M_{i,s} \in \PF_1(S)\}|$. Clearly, there are at most $\mu_1+\ldots+\mu_5$ elements $f \in \PF_1(S)$ such that a row of $\RF^+(f)$ contains $3$ zeroes. In order to have a bound for the number of pseudo-Frobenius numbers $f \in \PF_1(S)$ for which a row of $\RF^-(f)$ contains three zeroes, we do the following observations, where $h$ is the minimum for which $f_h \neq \F(S)$:
\begin{itemize}
\item  If $f_h - f = \lambda n_j - n_h$ for some index $j \neq h$ and $f \in \PF_1(S)$, then Proposition \ref{minimum index} implies $f_\ell - f = f_h-n_\ell+n_h - f = (f_h-f) - n_\ell + n_h=\lambda n_j - n_h-n_{\ell}+n_h=\lambda n_j -n_{\ell}$. Hence, we do not need to count when the $h$-th row of $\RF^-(f)$ has three zeroes because in this case also the $\ell$-th row has three zeroes (note also that if $j=\ell$ we get a contradiction). 
\item By Remark \ref{first zero}, for every $f \in \PF_1(S)$ the $(\ell,h)$-entry of $\RF^-(f)$ is zero, and so it is not possible that $\F(S)-f=\lambda n_h-n_\ell$ for some $\lambda \in \mathbb{N}$. 
\end{itemize}
Therefore, the elements $f \in \PF_1(S)$ for which a row of $\RF^-(f)$ contains three zeroes are at most $20-4-1=15$. Hence, we get $|\PF_1(S)|\leq \mu_1+\cdots+\mu_5+15$.

Now the result is clear when $\mu_1+\ldots+\mu_5 \leq 16$. The missing cases are when the values of $\mu_1, \ldots, \mu_5$ (up to order) are $(4,4,3,3,3)$, $(4,4,4,3,2)$, $(4,4,4,4,1)$, $(4,4,4,4,2)$, $(4,4,4,3,3)$, $(4,4,4,4,3)$ or $(4,4,4,4,4)$, but all these values are covered by the previous lemma.

Finally, if in the NG-vector of $S$ there are at least two entries different from $\F(S)$, then by Proposition \ref{minimum index} there is an element $f_{h'}$ in the NG-vector for which $\F(S)-f_{h'}=\lambda n_j - n_{h'}$ for some $j$ and $h' \neq h, \ell$. Since $f_{h'} \notin \PF_1(S)$, the elements $f \in \PF_1(S)$ for which a row of $\RF^-(f)$ contains three zeroes are at most $20-4-1-1=14$ and then $|\PF_1(S)|\leq \mu_1+\cdots+\mu_5+14$. Now it is enough to proceed as above. 
\end{proof}

\subsection{Bounds for the type}

As a consequence of Propositions \ref{PF2} and \ref{PF1} we can easily prove our main result.

\begin{theorem} \label{main}
If $S$ is a $5$-generated nearly Gorenstein numerical semigroup that is not almost symmetric, then $t(S) \le 40$. 
\end{theorem}

\begin{proof}
Recall that at least two entries of the NG-vector are equal by Corollary \ref{Distinct}. Moreover, $|\PF_2(S)|\leq 6$ by Proposition \ref{PF2}. 

If in the NG-vector of $S$ there are at least two entries different from $\F(S)$, then Proposition \ref{PF1} implies that
$$t(S) = |\{f_1,f_2,f_3,f_4,f_5\}| + |\PF_1(S)| + |\PF_2(S)| \le 4+30+6= 40.$$ 

On the other hand, if in the NG-vector of $S$ there is at most one entry different from $\F(S)$, then again Proposition \ref{PF1} implies that
\[t(S) = |\{f_1,f_2,f_3,f_4,f_5\}| + |\PF_1(S)| + |\PF_2(S)| \le 2+31+6= 39. \qedhere \]
\end{proof}

We end this section showing that the bound can be reduced in some particular case.

\begin{theorem} \label{3distinct}
Let $S$ be a $5$-generated nearly Gorenstein numerical semigroup with {\rm NG}-vector $(f_1, f_2, f_3, f_4,f_5)$, and assume that $f_1$, $f_2$, $f_3$ are pairwise distinct. Then, $t(S) \leq 5$.
\end{theorem}

\begin{proof}
Let $S$ be generated by $n_1 < n_2 < \dots < n_5$ and let $f \in \PF(S) \setminus \{f_1,f_2,f_3\}$. By Proposition \ref{f1-fj} we have $f_1-f=-n_1+a_4 n_4+a_5 n_5$ for some $a_4, a_5 \in \mathbb{N}$, $f_2-f_1=n_1-n_2$, and $f_3-f_1=n_1-n_3$. Therefore, since $f_i-f=(f_i-f_1)+(f_1-f)$ for $i=2,3$, we get
\begin{align*}
    f_1-f=-n_1+a_4 n_4+a_5 n_5, \\
    f_2-f=-n_2+a_4 n_4+a_5 n_5, \\
    f_3-f=-n_3+a_4 n_4+a_5 n_5.
\end{align*}
Therefore, these three equalities give the first three rows of an $\RF^-$--matrix for $f$. 
Clearly, at least one between $a_4$ and $a_5$ has to be positive. If $a_4 \neq 0$, then Lemma \ref{Coppie} implies that in the fourth row of every $\RF^+$--matrix for $f$ there are three zeroes, and more precisely $f=M_{4,5}=\lambda_{4,5}n_5-n_4$. Similarly, if $a_5 \neq 0$, then $f=M_{5,4}=\lambda_{5,4}n_4-n_5$. This means that $\PF(S) \subseteq \{f_1,f_2,f_3,M_{4,5},M_{5,4}\}$.
\end{proof}

\begin{example} Let $S=\langle 13, 45, 72, 79, 99 \rangle$. Using the {\rm NumericalSgps} package \cite{DGM} of {\rm GAP} \cite{GAP}, it is possible to see that $S$ is nearly Gorenstein and admits only two {\rm NG}-vectors: $(244, 212, 244, 244, 244)$ and $(244, 212, 185, 244, 244)$. Using the former we cannot apply Theorem \ref{3distinct}, but if we use the latter we get immediately that $t(S)\leq 5$. Indeed, in this case $\PF(S)=\{59, 185, 212, 244\}$. Note also that, according to the proof of the previous theorem, we have
\begin{align*}
    f_1-59=-n_1+2 n_5, \\
    f_2-59=-n_2+2 n_5, \\
    f_3-59=-n_3+2 n_5.
\end{align*}
and $59=2n_4 -n_5$.
\end{example}

\section{On the type of numerical semigroups}

Given a nearly Gorenstein numerical semigroup $S$ with embedding dimension $\nu$ and NG-vector $(f_1,f_2, \dots, f_\nu)$, if $S$ is not almost symmetric and $\nu=5$, in the previous section we proved that its type is at most $40$, by bounding the cardinalities of the two sets $\PF_1(S)$ and $\PF_2(S)$. 
In general, by definition we have that $|\PF_1(S)| \le 2\nu(\nu-1)$. Therefore, in order to figure out whether $t(S)$ is bounded by a function of $\nu$ we need to study $\PF_2(S)$. For $\nu=5$, the crucial fact in solving this case was that the elements of $\PF_2(S)$ have a specific form.
	
Generally speaking, information on the ``factorizations'' of a pseudo-Frobenius number $f \in \PF(S)$ (more formally, on the factorizations of $f+n_i$ for any $i$) is encoded by its RF-matrices $\RF^{+}(f)$ and $\RF^{-}(f)$; in particular, the distribution of zeroes on $\RF^+(f)$ gives information on how many generators appear in different factorizations of $f$. Moreover, since the number of possible distributions of zeroes of $\RF^+(f)$ and $\RF^-(f)$ is bounded for fixed $\nu$, then $t(S)$ is bounded by a function of $\nu$ if and only if for any admissible distribution of zeroes among the entries of two $\nu \times \nu$ matrices, there are at most a bounded number of $f \in \PF(S)$ with that distribution of zeroes among the entries of $\RF^+(f)$ and $\RF^-(f)$.
	
For $f \in \PF_2(S)$, any of the $2\nu$ rows of $\RF^+(f)$ and $\RF^-(f)$ must contain at most $\nu-3$ zeroes; on the other hand, by Lemma \ref{Coppie} there must be at least $\nu(\nu-1)$ zeroes among the entries of $\RF^+(f)$ and $\RF^-(f)$. 
	
For $\nu=5$ we have $2\nu(\nu-3)=20=\nu(\nu-1)$, and this forces specific distributions of zeroes for $\RF^+(f)$ and $\RF^-(f)$ (namely, each row and column having exactly $\nu-3=2$ zeroes), and further considerations on the nearly-Gorenstein property allowed us to reduce the number of admissible distributions to three, while Proposition \ref{PF2} implies that every such distribution can only be achieved by at most two elements of $\PF_2(S)$.
	
	On the other hand, for $\nu \ge 6$ we have $\nu(\nu-1) < 2\nu(\nu-3)$ and the difference between these two values is increasingly large, implying that there could be many more possible distributions of zeroes to be considered. Furthermore, while these distributions have different forms, it is not clear whether any given distribution can appear in row-factorization matrices of a bounded number of elements of $\PF(S)$.

In general, for numerical semigroups with fixed embedding dimension $\nu$ it is known that the type is not bounded by any function of $\nu$, and there are several such examples in literature. Many of them can be constructed by using RF-matrices; the key idea is to impose that there are an unbounded number (depending on a parameter) of pseudo-Frobenius numbers with RF-matrices having the same distribution of zeroes by using simple relations on the generators, and then deduce the generators from the non-zero coefficients of these RF-matrices (which will depend on the parameter as well).
We outline this idea with an example.
  
\begin{example}
		We want to build a family of numerical semigroups with embedding dimension $6$ and unbounded type. Let $S=\langle n_1,\ldots,n_6 \rangle$ be an element of this family; we assume that $n_4-n_1=n_5-n_2=n_6-n_3=d$, for some $d \in \mathbb{Z}^+$. Let $f \in \PF(S)$. Then $$f+\lambda d=f+\lambda(n_4-n_1)=f+\lambda(n_5-n_2)=f+\lambda(n_6-n_3) \text{   for } \lambda \in \mathbb{Z}^+.$$ We wish to leverage this equation to deduce (eventual) RF-matrices for $f+ \lambda d$, with $\lambda$ a positive integer, in function of a RF-matrix of $f$. In order to do so, we fix  $\RF^+(f)$, for instance by imposing that $$f=(T+1)n_2-n_1=Tn_3-n_2=(T+2)n_1-n_3,$$
		with $T \in \mathbb{Z}^+$ (this will be our parameter to construct the family). 
		These three relations actually lead to
		
		$$f=Tn_2+n_5-n_4=(T-1)n_3+n_6-n_5=(T+1)n_1+n_4-n_6,$$
		and hence
		$$\RF^+(f)=\begin{pmatrix}
			-1 & T+1 & 0 & 0 & 0 & 0 \\
			0 & -1 & T & 0 & 0 & 0  \\
			T+2 & 0 & -1 & 0 & 0 & 0 \\
			0 & T & 0 & -1 & 1 & 0 \\
			0 & 0 & T-1 & 0 & -1 & 1 \\
			T+1 & 0 & 0 & 1 & 0 & -1
		\end{pmatrix}.$$
	For $\lambda \in \mathbb{Z}^+$, we can deduce the first row of $\RF(f+\lambda d)$ from that of $\RF^+(f)$ by observing that $$f+\lambda d = ((T+1)n_2-n_1)+\lambda (n_5-n_2)=((T+1)-\lambda)n_2+\lambda n_5 - n_1.$$ Arguing in the same way, we can deduce that, if $f+\lambda d \not \in S$ and $\lambda \le T-1$, then $f+\lambda d$ will have a RF-matrix of the form 
		$$\RF^+(f + \lambda d)=\begin{pmatrix}
			-1 & T+1-\lambda & 0 & 0 & \lambda & 0 \\
			0 & -1 & T-\lambda & 0 & 0 & \lambda  \\
			T+2-\lambda & 0 & -1 & \lambda & 0 & 0 \\
			0 & T-\lambda & 0 & -1 & \lambda+1 & 0 \\
			0 & 0 & T-1-\lambda & 0 & -1 & \lambda+1 \\
			T+1-\lambda & 0 & 0 & \lambda+1 & 0 & -1
		\end{pmatrix},$$
		
		and thus $f+\lambda d \in \PF(S)$. Notice that all these $\RF^+(f+\lambda d)$ will have the same distribution of zeroes. 
		Next, we deduce the generators in function of $T$ and $d$. The equation $$f=(T+1)n_2-n_1=Tn_3-n_2=(T+2)n_1-n_3$$ can be solved in $n_1,n_2,n_3$, yielding for instance the solution 
		
		\begin{eqnarray*}
			n_1  & = & k[(T+1)^2+1] \\ 
			n_2 & = & k[(T+1)^2+T] \\
			n_3 &=&  k[(T+1)^2+2T+4] \\
			n_i & = & n_{i-3}+d \text{   for   } i = 4,5,6 \\
			f & = & k[T(T+1)(T+2)-1] 
		\end{eqnarray*}
		
		with $d,k,T \in \mathbb{Z}^+$.
		
		In order to obtain the desired family, we need to check that, for some suitable values of $d,k$, $S=\langle n_1,\ldots,n_6 \rangle$ is a numerical semigroup (that is, $\gcd(n_1,\ldots,n_6)=1$) and that $f+\lambda d \not \in S$ for $\lambda=0,\ldots,T-1$ (that is, $f+\lambda d$ cannot be written as $\alpha_1n_1+\cdots+\alpha_6n_6$, with $\alpha_1,\ldots,\alpha_6 \in \mathbb{N}$).
		
		For instance, it can be proved by elementary (although long) computations that this is the case if $k \ge T$, $d \ge T^2$, $\gcd(d,k) =1$ and $T \not \equiv 1 \pmod{5}$. 		
		Then this will yield a family of numerical semigroups $S_T$ with embedding dimension $6$ and such that $\{f,f+d,\ldots,f+(T-1)d\} \subseteq \PF(S_T)$.	
	\end{example}
	
	This idea can be used to recreate known examples from the literature.
	
	\begin{example} \label{Backelin}
		The classical example of Backelin \cite{FGH} $S_T=\langle s,s+3,s+3T+1,s+3T+2 \rangle$, for $T \ge 2$ and $s=(3T+2)^2+3$, can be obtained from the outlined method by assuming $n_4-n_3=1, n_2-n_1=3$, and $f=(3T+3)n_4-n_1$. Moreover, assume that for every $\lambda=1,\ldots,T$, $f-3\lambda \in \PF(S)$ and $$\RF^+(f-3\lambda)=\begin{pmatrix}
			-1 & 0 & 3\lambda & 3T+3-3\lambda \\
			0 & -1 & 3(\lambda-1) & 3T+3-3(\lambda-1) \\
			T+5-\lambda & 2T-1-\lambda & -1 & 0 \\
			2T+3+\lambda & T-\lambda & 1 & -1 \\
		\end{pmatrix}.$$
	\end{example}

	In our context, the issue lies in finding families of numerical semigroups which are nearly Gorenstein; namely, it seems difficult to compute a NG-vector from this information. The same problem arises even for almost-symmetric numerical semigroups. In fact, it is possible to construct numerical semigroups $S$ having an arithmetic progression of pseudo-Frobenius number, and that progression can be tweaked as to make sure that the sum of opposite terms is equal to $\F(S)$. However, this construction also gives no restriction on the other pseudo-Frobenius numbers, thus making difficult to verify whether $S$ is almost-symmetric or not. For instance, the following example was found by using this idea and checking computationally for almost symmetric numerical semigroups; however, it is not clear whether this specific example can be extended to an infinite family of semigroups.

\begin{example} \cite[Example 19]{M5} \label{Big AS}
Let $S=\langle 455,497,574,589,631,708 \rangle$ be a numerical semigroup. We have $$\PF(S)=\{ 3079, 3289, 3521, 3655, 3674, 3789, 3923, 4057, 4172, 4191, 4325, 4557, 4767, 7846 
	\}.$$ Then $S$ is an almost symmetric numerical semigroup with type $t(S) = 14 > 2\nu$. Further, $\PF(S)$ contains the arithmetic progression $\{3521,3655,3789,3923,4057,4191,4325\}$ of ratio $134$, and
	 $$\RF^+(3521+134\lambda) = \begin{pmatrix}
		-1 & 8-\lambda & 0 & 0 & \lambda & 0 \\
		0 & -1 & 7-\lambda & 0 & 0 & \lambda  \\
		9-\lambda & 0 & -1 & \lambda & 0 & 0 \\
		0 & 7-\lambda & 0 & -1 & \lambda+1 & 0 \\
		0 & 0 & 6-\lambda & 0 & -1 & \lambda+1 \\
		8-\lambda & 0 & 0 & \lambda+1 & 0 & -1
	\end{pmatrix}
	$$
	
	for $\lambda=1,\ldots,6$. Thus, the pseudo-Frobenius numbers $\{3655,3789,3923,4057,4191\} \subseteq \PF(S)$ (obtained for $\lambda=1,\ldots,5$) have the same distribution of zeroes on their $\RF$-matrices.
\end{example}

In \cite{Ba} Barucci asked if an almost symmetric numerical semigroup for which $t(S) \ge 2 \nu$ exists, and Example \ref{Big AS} is the first example of this kind found in literature.

\begin{remark} \rm \label{Numerical duplication}
We notice also that Example \ref{Big AS} can be used to show that the type of almost symmetric numerical semigroups with embedding dimension $\nu(S)$ cannot be bounded by $2\nu(S)+a$ with $a \in \mathbb{Z}$. To explain why this is true, we need to use the notion of numerical duplication introduced in \cite{DS}. Given a numerical semigroup $S$, an ideal $E$ of $S$ and an odd integer $b \in S$, the {\rm numerical duplication} of $S$ with respect to $E$ and $b$ is the numerical semigroup
\[
S\!\Join^b\!E=2\cdot S \cup (2 \cdot E+b),
\]
where $2\cdot X=\{2x \mid x \in X\}$. If $S$ is almost symmetric and $E=M(S):=S\setminus \{0\}$ is the maximal ideal of $S$, then by \cite[Theorem 3.4]{BDS} and \cite[Proposition 2.9]{BDS2} follow that also $S\!\Join^b\!M(S)$ is almost symmetric with type $t(S\!\Join^b\!M(S))=2t(S)+1$ for every odd integer $b \in S$. Moreover, $\nu(S\!\Join^b\!M(S))=2\nu(S)$ by \cite[Proposition 2.3]{BDS}, where $\nu(\cdot)$ denotes the embedding dimension of a numerical semigroup. Thus, if $S$ is almost symmetric, then $S_1:=S\!\Join^b\!M(S)$ is an almost symmetric numerical semigroup with $t(S_1)-2\nu(S_1)=2t(S)+1-4\nu(S) =2(t(S)-2\nu(S))+1$. More generally, if for every $i$ we define $S_i:=S_{i-1}\!\Join^{b_{i-1}}\!M(S_{i-1})$ with $b_{i-1} \in S_{i-1}$ odd, then $S_i$ is almost symmetric and it is easy to prove by induction that $t(S_i)-2\nu(S_i) =2^i (t(S)-2\nu(S))+2^i-1$. When $S$ is the semigroup in Example \ref{Big AS}, $t(S)-2\nu(S)=2$ is positive, and then for every $a \in \mathbb{Z}$ we can construct an almost symmetric numerical semigroup $S_i$ for which $t(S_i)-2\nu(S_i)>a$, i.e., $t(S_i)>2\nu(S_i)+a$.
\end{remark}

\section*{Acknowledgments}
The authors would like to thank the anonymous referee for some helpful comments that improved the quality of the manuscript.


\begin{thebibliography}{Dillo 83}
\footnotesize

\bibitem{ADG} A. Assi, M. D'Anna, P.A. Garc\'ia-S\'anchez, {\em Numerical semigroups and applications}, RSME Springer Series, Vol. 3, 2020.


\bibitem{Ba} V. Barucci, {\em Almost symmetic property in semigoups and rings}, presentation at IMNS 2016, Levico Terme (Italy), available at \href{https://www.ugr.es/~imns2010/2016/Slides/barucci-imns2016.pdf}{https://www.ugr.es/$\sim$imns2010/2016/Slides/barucci-imns2016.pdf}.

\bibitem{BDS} V. Barucci, M. D'Anna, F. Strazzanti, {\em A family of quotients of the Rees Algebra}, Comm. Algebra {\bf 43} (2015), no. 1, 130-142.

\bibitem{BDS2} V. Barucci, M. D'Anna, F. Strazzanti, {\em Families of Gorenstein and almost Gorenstein rings}, Ark. Mat. {\bf 54} (2016), no. 2, 321--338.


\bibitem{BF} V. Barucci, R. Fr\"oberg, {\em One-dimensional almost Gorenstein rings},  J. Algebra {\bf 188} (1997), 418--442.



\bibitem{DS} M. D'Anna, F. Strazzanti, {\em The numerical duplication of a numerical semigroup}, Semigroup Forum {\bf 87} (2013), no. 1, 149-160.

\bibitem{DGM} M. Delgado, P.A. Garc\'ia-S\'anchez, J. Morais, {\em NumericalSgps, a package for numerical semigroups}, GAP package, Version 1.3.1 (2022). 

\bibitem{D} S. Ding, {\em A note on the index of Cohen-Macaulay local rings}, Comm. Algebra {\bf 21} (1993), no. 1, 53--71.

\bibitem{E} K. Eto, {\em Almost Gorenstein monomial curves in affine four space}, J. Algebra {\bf 488} (2017), 362--387.

\bibitem{FGH} R. Fr\"oberg, C. Gottlieb, R. H\"aggkvist, {\em On numerical semigroups}, Semigroup Forum {\bf 35} (1986), 63--83.

\bibitem{GAP} The GAP Group, {\em GAP -- Groups, Algorithms, and Programming}, Version 4.10.2 (2019). 


\bibitem{H} J. Herzog, {\em Generators and relations of abelian semigroups and semigroup rings}, Manuscr. Math. {\bf 3} (1970), 175--193.

\bibitem{HHS} J. Herzog, T. Hibi, D.I. Stamate, {\em The trace of the canonical module}, Isr. J. Math. {\bf 233} (2019), 133--165.

\bibitem{HHS2} J. Herzog, T. Hibi, D.I. Stamate, {\em Canonical trace ideal and residue for numerical semigroup rings}, Semigroup Forum {\bf 103} (2021), 550--566.

\bibitem{HW} J. Herzog, K.-i. Watanabe, {\em Almost symmetric numerical semigroups}, Semigroup Forum {\bf 98} (2019), no. 3, 589--630.

\bibitem{HV} C. Huneke, A. Vraciu, {\em Rings that are almost Gorenstein}, Pacific J. Math. {\bf 225} (2006), 85--102.


\bibitem{M5} A. Moscariello, {\em On the boundedness of the type of an almost Gorenstein monomial curve in $\mathbb{A}^5$}, Comm. Algebra {\bf 51} (2023), no. 3, 1179--1185.

\bibitem{M} A. Moscariello, {\em On the type of an almost Gorenstein monomial curve}, J. Algebra {\bf 456} (2016), 266--277.

\bibitem{MS} A. Moscariello, F. Strazzanti, {\em Nearly Gorenstein vs almost Gorenstein affine monomial curves}, Mediterr. J. Math. {\bf 18} (2021), article number: 127.

\bibitem{N} H. Nari, {\em Symmetries on almost symmetric numerical semigroups}, Semigroup Forum {\bf 86} (2013), no. 1, 140--154.

\bibitem{Nu} T. Numata, {\em Almost symmetric numerical semigroups generated by four elements}, Proceedings of the Institute of Natural Sciences, Nihon University, {\bf 48} (2013), 197--207.


\bibitem{RG} J.C. Rosales, P.A. Garc\'ia-S\'anchez, {\em Numerical Semigroups}, Springer Developments in Mathematics, Vol 20, 2009.


\bibitem{S} D.I. Stamate, {\em Betti numbers for numerical semigroup rings}, in Multigraded Algebra and Applications, NSA 2016, Springer Proceedings in Mathematics \& Statistics {\bf 238} (2018).

\end{thebibliography}
\end{document}